\documentclass{article}

\usepackage[cmtip,arrow]{xy}
\usepackage{amsmath,amssymb,enumerate,pb-diagram,pb-xy}

\usepackage[applemac]{inputenc}

\def\A{{\mathbb A}}
\def\adm{{\rm adm}}
\def\al{\alpha}
\def\Ann{\operatorname{Ann}}
\def\C{{\mathbb C}}
\def\CB{{\cal B}}
\def\CH{{\cal H}}
\def\df{\ \stackrel{\mbox{\rm\tiny def}}{=}\ }
\def\e{\emph}
\def\End{\operatorname{End}}
\def\eps{\varepsilon}

\def\ga{\gamma}

\def\Hom{\operatorname{Hom}}
\def\Id{{\rm Id}}

\def\la{\lambda}

\def\ph{\varphi}
\def\prf{{\bf Proof: }}

\def\qed{\phantom{}\ensuremath{\hfill\square}}

\def\sm{\smallsetminus}
\def\tr{\operatorname{tr}}
\def\vol{{\rm vol}}
\def\what{\widehat}

\def\({\left(}
\def\){\right)}
\def\={{\ =\ }}

\newcommand{\tto}[1]{\stackrel{#1}{\longrightarrow}}

\newcommand{\norm}[1]{|\hspace{-1pt}| #1|\hspace{-1pt}|}
\renewcommand{\sp}[1]{\left\langle #1\right\rangle}
\newcommand{\ol}[1]{\overline{#1}}

\newtheorem{theorem}{Theorem}[section]

\newtheorem{lemma}[theorem]{Lemma}

\newtheorem{proposition}[theorem]{Proposition}

\newtheorem{defi}[theorem]{Definition}

\begin{document}

\pagestyle{myheadings} \markright{TENSOR PRODUCT THEOREM}

\title{A general tensor product theorem}
\author{Anton Deitmar}
\date{}
\maketitle

{\bf Abstract:} We show that every admissible irreducible representation of a product of two locally compact groups is a tensor product of admissible irreducible representations of the factors.

\section*{Introduction}
In the classical case the tensor product theorem is formulated for either totally disconnected groups or Lie groups and the proof by D. Flath \cite{Flath}, which is found in the books until this day, uses different methods for each case.
The tensor product theorem is used in the theory of automorphic forms, to prove that every irreducible admissible representation of an adele group $G(\A)$ is an infinite product of local representations.
In this note we give a unified proof which is simpler and applies to all locally compact groups.
The proof essentially contains no new ideas. We only streamline and unify the existing proofs by using a more universal concept of Hecke algebras.
This proof will also appear in the upcoming book by the author on Automorphic Forms.

\section{Preliminaries}

For a locally compact group $G$ we denote by $\hat G$ the \e{unitary dual}, i.e., the set of isomorphy classes of irreducible unitary representations of $G$.

For a representation $(\pi,V_\pi)$ of the compact group  $K$ and an irreducible representation $(\tau,V_\tau)$ of $K$ 
we define the \e{$\tau$-isotype} as the sum of all subrepresentations of $\pi$, which are isomorphic to $\tau$.
If $\pi$ is unitary, then $V_\pi$ is the direct sum of its isotpyes
$
V_\pi\=\bigoplus_{\tau\in\hat K}V_\pi(\tau).
$

For a given representation $\tau\in\hat K$ let
$$
e_\tau(k)\=(\dim \tau)\ \ol{\tr(\tau(k))},\qquad k\in K.
$$

\begin{lemma}\label{7.5.10}
The function $e_\tau$ is an idempotent in the convolution algebra $C(K)$, i.e., one has $e_\tau *e_\tau=e_\tau$.
Let $(\pi,V_\pi)$ be a representation of $K$. Then the map
$$
P_\tau\=\int_Ke_\tau(k)\pi(k)\,dk
$$
is a projection onto the isotype $V_\pi(\tau)$.
Here we have normalized the Haar measure on $K$ so that $\vol(K)=1$.
If $\pi$ is unitary, then $P_\tau$ is the orthogonal projection onto $V_\pi(\tau)$.

If $\ga$ is another irreducible unitary representation of $K$ which is not isomorphic to $\tau$, then 
$
e_\tau*e_\ga\= 0.
$
\end{lemma}

\prf
This is a direct consequence of the Peter-Weyl Theorem \cite{HA2}.
\qed

Let $F\subset\hat K$ be a finite subset.
Then the function
$$
e_F\df\sum_{\tau\in F}e_\tau
$$
is an idempotent in $C(K)$, because
$$
e_F*e_F\=\sum_{\tau,\ga\in F}e_\tau*e_\ga\=\sum_{\tau\in F} e_\tau\= e_F.
$$

\section{Main Theorem}
Let $G$ be a locally-compact group and $K\subset G$ a compact subgroup.
A unitary representation $(\pi,V_\pi)$ of $G$ is called \e{$K$-admissible}, if the representation $\pi|_K$ decomposes with finite multiplicities, i.e., if for every $\tau\in\what K$ the multiplicity
$$
[\pi:\tau]\=\dim\Hom_K(V_\tau,V_\pi)\=\dim\Hom_K(V_\pi,V_\tau)
$$
is finite.
The representation $\pi$ is $K$-admissible if and only if every $K$-isotype $V_\pi(\tau)$ is finite-dimensional.
We denote by \e{$\what G_K$} the set of all isomorphy classes of irreducible unitary $G$-representations which are $K$-admissible.
A representation $\pi$ of $G$ is called \e{admissible}, if there is a compact subgroup $K$ such that $\pi$ is $K$-admissible. 
We write $\what G_\adm$ for the set of isomorphy classes of irreducible unitary admissible representations of $G$.

\begin{theorem}[Main Theorem]\label{4.5.9}
Let $G$ and $H$ be locally compact groups.
For any two $\pi\in\what G_\adm$ and $\eta\in\what H_\adm$ the tensor product $\pi\otimes\eta$ is in $\what{G\times H}_\adm$ and the tensor map $(\pi,\eta)\mapsto \pi\otimes\eta$ is a bijection
$$
\what G_\adm\times\what H_\adm\ \tto\cong\ \what{G\times H}_\adm.
$$
More precisely, given two compact subgroups $K\subset G$ and $L\subset H$,
the tensor map map is a bijection  $\what G_K\times\what H_L\to\what{G\times H}_{K\times L}$.
\end{theorem}

The injectivity of the tensor map follows from the next lemma.

\begin{lemma}\label{7.6.3}
Let $G$ be a locally compact group and let $(\eta,V_\eta)$ be an irriducible unitary representation.
Let $W$ be a Hilbert space. 
Then every closed $G$-stable subspace of $V_\eta\otimes W$ is of the form $V_\eta\otimes W_1$ for a closed subspace $W_1$ of $W$.
Here the tensor product refers to the tensor product in the category of Hilbert spaces, i.e., the Hilbert-space completion of the algebraic tensor product.

In particular, it follows that for two locally compact groups  $G,H$ and irreducible unitary representations $\pi,\tau$ of $G$ and $H$ respectively, the representation $\pi\otimes\tau$ of $G\times H$ is irreducible and that $\pi\otimes\tau\ \cong\ \pi'\otimes\tau'$ implies  $\pi\cong\pi'$ and $\tau\cong\tau'$.
\end{lemma}

\prf
Let $U\subset V_\eta\otimes W$ be a closed, $G$-stable subspace.
Then the orthogonal projection $P$ with image $U$ is a bounded operator on $V_\eta\otimes W$ which commutes with all  $\eta(g)$, $g\in G$.
We show that every such operator $T$ is of the form $1\otimes S$ for an operator $S\in\CB(W)$.
So let $T$ be a bounded operator on $V_\eta\otimes W$ with 
$$
T(\eta(g)\otimes 1)\=(\eta(g)\otimes 1)T
$$
for every $g\in G$.
Let $(e_i)_{i\in I}$ be an orthonormal basis of $W$ and for $j\in I$ let $P_j$ be the orthogonal projection onto the subspace $\C e_j$.
For $i,j\in I$ consider the operator $T_{i,j}$ given as composition:
$$
V_\eta\to V_\eta\otimes e_i\tto T V_\eta\otimes W\tto{1\otimes P_j} V_\eta\otimes e_j\to V_\eta.
$$
This operator is bounded and commutes with the $G$-action.
By the Lemma of Schur there exists a complex number 
$a_{i,j}$ with $T_{i,j}=a_{i,j}\Id$.
Since $\sum_jP_j=\Id_W$, we get
$$
T(v\otimes e_i)\=\sum_ja_{i,j}v\otimes e_j\=v\otimes S(e_i),
$$
where $S(e_i)=\sum_ja_{i,j}e_j$ and the sum converges in $W$.
As $T$ is continuous, $S$ extends to a unique continuous operator $S:W\to W$, satisfying
$$
T(v\otimes w)\=v\otimes S(w)
$$
for all $v\in V_\eta$, $w\in W$.
This means $T=1\otimes S$, as claimed.

We can apply this to the orthogonal projection $P$ onto the invariant subspace $U$, which therefore is of the form $P=1\otimes P_1$.
Then $P_1$ is again an orthogonal projection.
Let $W_1$ be the image of $P_1$, then $U= V_\eta\otimes W_1$.

We next show the irreducibility of $\pi\otimes\tau$.
The first part implies that a closed, $G\times H$-stable subspace $U$ has to be of the form $U=V_\pi\otimes U_\tau$ as well as of the form $U=U_\pi\otimes V_\tau$ with closed subspaces $U_\pi\subset V_\pi$ and $U_\tau\subset V_\tau$.
These subspaces have to be invariant themselves, hence the full spaces.

Finally let $\pi\otimes\tau\cong\pi'\otimes\tau'$, 
so there is a unitary isomorphism $T:V_\pi\otimes V_\tau\to V_{\pi'}\otimes V_{\tau'}$ commuting with the $G\times H$-operation.
Choose $w\in V_\tau$ with $\norm w=1$.
Then $V_\pi\otimes w$ is an irreducible $G$-subrepresentation, therefore there is $w'\in V_{\tau'}$ mit $\norm{w'}=1$, such that $T(V_\pi\otimes w)=V_{\pi'}\otimes w'$.
The map
$$
V_\pi\to V_\pi\otimes w\tto T V_{\pi'}\otimes w'\to V_{\pi'}
$$ 
is a unitary $G$-isomorphism, so we get $\pi\cong\pi'$ and analogously $\tau\cong\tau'$.
\qed

\section{The Hecke algebra}
Let $G$ be a locally compact group.
The set $C_c(G)$ of all continuous functions of compact support is a complex algebra under the convolution product
$$
f*g(x)\=\int_Gf(y)g(y^{-1}x)\,dy\=\int_Gf(xy)g(y^{-1})\,dy,\quad f,g\in C_c(G).
$$

Let $K\subset G$ be a compact subgroup.
For $\al,\beta\in C(K)$, there is the $K$-convolution
$$
\al *\beta(k)\=\int_K\al(l)\beta(l^{-1}k)\,dl
$$
making $C(K)$ a $\C$-algebra as well.
We can also define a convolution between functions on $G$ and $K$ as follows.
For $\al\in C(K)$ and $f\in C_c(G)$ we define
$$
\al*f(x)\=\int_K\al(k)f(k^{-1}x)\,dk\quad \text{and}\quad f*\al(x)\=\int_Kf(xk)\al(k^{-1})\,dk.
$$
We normalize the Haar-measure of the compact group $K$ so that $\vol(K)=1$.
If $K$ is not a nullset in $G$, we insist that the Haar measure of $G$ be normalized also to give $K$ the volume $1$.
Hence in this case the Haar measures of $G$ and $K$ coincide on $K$. 

\begin{lemma}
With this normalization the convolution products are compatible in the sense that 
$$
f*(g*h)\=(f*g)*h.
$$
holds for all combinations of $f,g,h$ being functions in either  $C_c(G)$ or $C(K)$.
Also we always have
$
(f*g)^*\=g^**f^*,
$
where $f^*(x)=\Delta(x^{-1})\ol{f(x^{-1})}$ and $\Delta(x)$ is the modular function of $G$.
\end{lemma}

\prf
A computation.\qed

Let $G$ be a locally compact group with a compact subgroup $K$.
A function $f$ on $G$ is called \e{$K$-finite}, 
if the set of all functions  $x\mapsto f(k_1xk_2)$, $k_1,k_2\in K$ spans a finite-dimensional vector space.
It is a simple observation that the convolution product of two $K$-finite functions is $K$-finite.
We define the \e{Hecke-algebra} $\CH=\CH_{G,K}$ 
of the pair $(G,K)$ as the convolution algebra of all $K$-finite functions in $C_c(G)$.

\begin{lemma}\label{Lem7.5.25}
Let $G$ be alocally compact group and $K$ a compact subgroup.
\begin{enumerate}[\rm (a)]
\item Let $I_K$ be the set of all finite subsets of $\hat K$.
For $F\in I_K$ set $e_F=\sum_{\tau\in F}e_\tau$ and
$$
C_F\df e_F* C_c(G)*e_F\=e_F* \CH *e_F.\index{$C_F$}
$$
Then $e_F$ is an idempotent in $C(K)$ and $C_F$ is a subalgebra of $\CH$.
The Hecke-algebra $\CH$ is the union of all these subalgebras.

If $F\subset F'$ in $I_K$, then $C_F\subset C_{F'}$.
If $F=\{\tau\}$ we also write $C_F=C_\tau$.
\item For a unitary representation $(\pi,V_\pi)$ of $G$ the space $\pi(\CH)V_\pi$ is dense in $V_\pi$.
\item The Hecke-algebra is a *-Algebra with involution $f^*(x)=\Delta(x^{-1})\ol{f(x^{-1})}$.
If $\pi$ is a unitary representation of $G$, then $\pi$
defines a *-representation of $\CH$.
For two unitary representations $\pi, \pi'$ of $G$ we have
$$
\pi\cong_G\pi'\ \Leftrightarrow\ \pi\cong_\CH\pi',
$$
so, $\pi$ and $\pi'$ are unitarily equivalent if and only if they define isomorphic $\CH$-modules.
\end{enumerate}
\end{lemma}

\prf
Part (a) is clear. We show (b).
For $h\in C(K)$ we write
$$
\pi(h)\=\int_K h(k)\pi(k)\,dk.
$$
Let $F\in I_K$ and $P_F=\pi(e_F)$.
Then $P_F^2=\pi(e_F)\pi(e_F)\=\pi(e_F*e_F)\=\pi(e_F)=P_F$.
So $P_F$ is a projection with image $V_\pi(F)=\bigoplus_{\tau\in F}V_\pi(\tau)$, which is the  \e{$F$-isotype} of $\pi$ and the kernel is $\bigoplus_{\tau\notin F}V_\pi(\tau)$.
Hence the union of all $V_\pi(F)$ with $F\in I_K$ is dense in $V_\pi$.
So it suffices to show that $\pi(C_F)V_\pi$ is dense in $V_\pi(F)$.
Let $v\in V_\pi(F)$ and let $\eps>0$.
Since $\pi(e_F)$ is continuous, there is $C>0$ such that $\norm{\pi(e_F)w}\le C\norm w$ for every $w\in V_\pi$.
Since the map $G\times V_\pi\to V_\pi$; $(g,v)\mapsto \pi(g)v$ is continuous, there is a neighborhoor $U$ of the unit in $G$, such that $x\in U\ \Rightarrow\ \norm{\pi(x)v-v}<\eps/C$.
Let $f\in C_c(G)$ with support in $U$ such that $f\ge 0$ and $\int_G f(x)\, dx\= 1$.
Then
$$
\norm{\pi(f)v-v}\=\norm{\int_Gf(x)(\pi(x)v-v)\,dx}
\le {\int_Gf(x)\norm{\pi(x)v-v}\,dx} <\eps/C.
$$
One has $e_F*f*e_F\in C_F\subset\CH$ and
$$
\norm{\pi(e_F*f*e_F)v-v}\=\norm{\pi(e_F)(\pi(f)v-v)}< \eps.
$$
This shows (b).
We continue with (c).
The closedness under * is clear.
Let $\pi$ be a unitary $G$-representation, then for $f\in\CH$ one has
\begin{eqnarray*}
\pi(f^*)&=&\int_G\Delta(x^{-1})\ol{f(x^{-1})}\pi(x)\,dx
\=\int_G\ol{f(x)}\pi(x^{-1})\,dx\\
&=&\int_G\ol{f(x)}\pi(x)^*\,dx\=\(\int_G{f(x)}\pi(x)\,dx\)^*\=\pi(f)^*.
\end{eqnarray*}
If $\pi$ and $\pi'$ are isomorphic as $G$-representations, then so they are as $\CH$-modules.
For the converse let $T: V_\pi\to V_{\pi'}$ be a unitary
$\CH$-isomorphism, so $T\pi(f)=\pi'(f)T$ for every $f\in\CH$.
We first remark, that the same must hold for $f\in C_c(G)$.
For this let $S=T\pi(f)-\pi'(f)T$, then for every finite subset $F$ of $\what K$,
$$
\pi'(e_F)S\pi(e_F)\=T\pi(\underbrace{e_F fe_F}_{\in\CH})-\pi'(e_Ffe_F)T\=0.
$$
Therefore $Sv=0$ for every vector $v\in V_\pi(F)$ and since the $V_\pi(F)$ span a dense subspace of $V_\pi$, we get $S=0$, so $T\pi(f)=\pi'(f)T$ for every $f\in C_c(G)$.

Let $\eps>0$ and let $x\in G$, $v\in V$. 
Then, by the continuity of the representations $\pi$ and $\pi'$ there exists a neighborhood $U$ of $x\in G$ such that $\norm{\pi(u)v-\pi(x)v}<\eps/2$ and $\norm{\pi'(u)Tv-\pi'(x)Tv}<\eps/2$.
Let $f\in C_c(G)$ have support inside $U$ and suppose $f\ge 0$ and $\int_Gf(x)\,dx=1$.
Since $T$ is unitary, we have
\begin{multline*}
\norm{T\pi(f)v-T\pi(x)v}
\=\norm{\pi(f)v-\pi(x)v}\\
\=\norm{\int_Gf(u)(\pi(u)v-\pi(x)v)\,du}\ \le\ \int_Gf(u)\norm{\pi(u)v-\pi(x)v}\,du\ <\ \eps/2
\end{multline*}
and likewise $\norm{T\pi(f)v-\pi'(x)Tv}=
\norm{\pi'(f)Tv-\pi'(x)Tv}<\eps/2$.
This implies $\norm{T\pi(x)v-\pi'(x)Tv}<\eps$.
Since $\eps>0$ and $v\in V_\pi$ are arbitrary, we conclude $T\pi(x)=\pi'(x)T$.
\qed

If $(\pi,V_\pi)$ is an irreducible unitary representation of $G$, then every $f\in C_c(G)$ defines a continuous linear operator $\pi(f)$ on  $V_\pi$.
For $0\ne v\in V_\pi$ the space  $\pi(C_c(G))v$ is a $G$-stable subspace of $V_\pi$.
As $\pi$ is irreducible, the space $\pi(\CH)v$ is a dense subspace.
 
\begin{proposition}\label{4.5.11}
Let $G$ be a locally compact group and $K\subset G$ a compact subgroup.
\begin{enumerate}[\rm (a)]
\item 
Let $(\pi,V_\pi)$ be an irreducible representation of $G$.
Let $F$ be a finite subset of $\hat K$.
Then the $F$-isotype $F$-Isotyp $V_\pi(F)=\bigoplus_{\tau\in F}V_\pi(\tau)$ is an irreducible $C_F$-module.
In particular, $\pi(\CH)V_\pi$ is an irreducible $\CH$-module.
\item 
Let $(\eta,V_\eta)$ be a unitary representation of $G$ and $F$ a finite subset of $\hat K$.
If $M$ is a finite-dimensional irreducible $C_F$-submodule of $V_\eta(F)$, then the $G$-subrepresentation of $\pi$ generated by $M$ is irreducible.
\item
Let $\eta,\pi$ be irreducible unitary representations of $G$, which are $K$-admissible. Then
$$
\eta\cong\pi\ \Leftrightarrow\ V_\eta(F)\ \cong\ V_{\pi}(F) \text{ for every } F\in I_K,
$$
where on the right we consider isomorphy as $C_F$-modules.
\end{enumerate}
\end{proposition} 

\prf
(a) To show irreducibility of $V_\pi(F)$, let $0\ne U\subset V_\pi(F)$ be a closed submodule.
Since $V_\pi$ is irreducible, the space $\pi(C_c(G))U$ is dense in $V_\pi$.
Let $P_F:V_\pi\to V_\pi(F)$ be the isotypical projection. 
For $h\in C(K)$ we write
$$
\pi(h)\=\int_Kh(k)\pi(k)\,dk.
$$
Then we have $P_F=\pi(e_F)$.
Let $f\in C_c(G)$. Then
$$
\pi(f)\pi(e_\tau)=\pi(f*e_\tau)\quad\text{und}\quad \pi(e_\tau)\pi(f)=\pi(e_\tau *f),
$$
as a calculation shows.

The projection $P_F$ is continuous, so the space $P_F(C_c(G) U)$ is dense in $V_\pi(F)$, therefore $\ol{P_F(C_c(G)U)}=V_\pi(F)$.
Since $P_F=\pi(e_F)$ we get
\begin{eqnarray*}
V_\pi(F)&=&\ol{P_F(C_c(G)U)}\=\ol{\pi(e_F)C_c(G)U}\\
&=&\ol{\pi(e_F*C_c(G)*e_F)U}\=\ol{\pi(C_F)U}\=\ol{U}\= U,
\end{eqnarray*}
which shows that $V_\pi(F)$ is irreducible.

We now show (b).
Let $m_0\in M\sm\{ 0\}$ and let
$
\Ann(m_0)\=\{a\in C_F:am_0=0\}
$
be its annihilator.
Then $\Ann(m_0)$ is a left ideal of $C_F$ and the map $a\mapsto am_0$ is a module isomorphism 
$$
C_F/\Ann(m_0)\to C_Fm_0.
$$
Since $M$ is finite-dimensional  and irreducible it follows $M=C_Fm_0$, so $M\cong C_F/J$ with $J=\Ann(m_0)$.
Let $U$ be the $G$-stable closed subspace of $V_\pi$ generated by $M$.
We show $P_F(U)=M$, where $P_F=\pi(e_F)$ is the orthogonal projection onto $M$.

\emph{Claim:} Let $\ol J$ be the annihilator $\Ann_{C_c(G)}(m_0)$ of $m_0$ in $C_c(G)$.
Then
$$
\ol J\=\ol J e_F \oplus \Ann_{C_c(G)}(e_F),
$$
where $\Ann_{C_c(G)}(e_F)$ is the annihilator of $e_F$ in $C_c(G)$, which is the set of all $f\in C_c(G)$ with $f*e_F=0$.

Ww show tha claim.
It is clear that $\ol Je_F\subset\Ann_{C_c(G)}(v_0)=\ol J$. Further $v_0=e_F v_0$ and so $\Ann_{C_c(G)}(e_F)\subset\Ann_{C_c(G)}(v_0)=\ol J$.
It remains to show that $\ol J$ lies in the right hand side.
Since $e_F$ is an idempotent, we have $C_c(G)=C_c(G)e_F\oplus \Ann_{C_c(G)}(e_F)$, since every $f\in C_c(G)$ can be written as $f=fe_F+(f-fe_F)$ and $f-fe_F$ lies in the annihilator of $e_F$.
Further 
$
\Ann_{C_c(G)}(e_F)\ \subset\ \Ann_{C_c(G)}(v_0)=\ol J
$, which shows the claim.

So we have
$$
C_c(G)v_0\cong C_c(G)/\ol J\cong C_c(G)e_F/ (\ol Je_F).
$$
This implies
$$
P_F(C_c(G)v_0)
\cong e_F C_c(G)e_F/(\ol Je_F)
\cong C_F/J\cong M.
$$
Therefore $P_F(U)=M$.
Let now $U'$ be a subrepresentation of $U$. Then $P_F(U')=0$ or $M$. In the first case $M\subset (U')^\perp$ and so $U\subset (U')^\perp$, since the space $(U')^\perp$ is a subrepresentation. 
This gives $U'=0$.
In the second case $M\subset U'$ and so $U'=U$.
Indeed we conclude that $U$ is irreducible.

(c) Let $(\eta,V_\eta),(\pi,V_\pi)$ be isomorphic representations, then so are the $C_F$-modules $V_\eta(F)$ and $V_\pi(F)$.
For the converse let for every $F$ be given a $C_F$-isomorphism 
$\phi_F:V_\eta(F)\to V_\pi(F)$.
According to the Lemma of Schur two isomorphisms for a given $F$ only differ by a multiple scalar.
So we can normalize all these isomorphisms in a manner that they extend each other, i.e., that we have 
$\phi_F=\phi_{F'}|_{V_\eta(F)}$, if $F\subset F'$.
This means that the  $\phi_F$ can be put together to a $\CH$-isomorphism
$$
\phi:\eta(\CH)V_\eta\tto\cong\pi(\CH)V_\pi.
$$
Both sides are dense sub vector spaces.
If we can show that $\phi$ is isometric, this map extends to an isomorphism 
of $G$-representations by Lemma \ref{Lem7.5.25}.
Let $0\ne v_0\in \eta(\CH)V_\eta$ and let $w_0=\phi(v_0)$.
We can assume that $\norm{v_0}=\norm{w_0}=1$ and we claim that $\phi$ is isometric.
For this we transport both inner  products to the same sode via $\phi$, so we have two inner products $\sp{.,.}_1$ and $\sp{.,.}_2$ on, say $\pi(\CH)V_\pi$ such that $v_0$ has norm $1$ in both inner products and the operation of $\CH$ is a $*$-operation in both inner products.
On the finite dimensional space $\pi(C_F) V_\pi$ these inner products must coincide by the Lemma of Schur.
This holds for every $F\in I_K$, so that $\phi$ is indeed isometric.
\qed

Let now $G,H$ locally compact with compact subgroups  $K\subset G$ and $L\subset H$.
Let $E$ be a finite subset of $\hat K$ and $F$ a finite subset of $\hat L$.
There is a natural Homomorphism
$$
\psi:C_c(G)\otimes C_c(H)\to G_c(G\times H)
$$
given by
$$
\psi(f\otimes g)(x,y)\=f(x)g(y).
$$
This induces a homomorphism
$$
C_E\otimes C_F
\to C_{E\times F}.
$$
This map is not surjective in general.
But this fact causes no harm, as we have

\begin{lemma}\label{Lem7.5.22}
If $M$ is a finite-dimensional irreducible $C_{E\times F}$-*-submodule of a unitary $G$-representation, then it is irreducible as $C_E\otimes C_F$-module.
\end{lemma}

\prf
We equip $C_{E\times F}\subset L^1(G\times H)$ with the topology of the $L^1$-norm.
Then $C_E\otimes C_F$ is dense in $C_{E\times F}$, as $C_c(G)$ is dense in $L^1(G)$.
The representation $\rho:C_{E\times F}\to\End_\C(M)$ is continuous, since $M$ comes from a unitary representation of $G\times H$.
So the image of $C_E\otimes C_F$ is dense in the (finite-dimensional) image of $C_{E\times F}$ in $\End(M)$, 
so these two images agree.
\qed

In order to show Theorem \ref{4.5.9} we only need the second assertion of the following lemma.

\begin{lemma}\label{Lem7.5.23}
\begin{enumerate}[\rm (a)]
\item Let $A$ be a  $\C$-algebra with unit and  $M$ a simple $A$-module, which is finite dimensional over $\C$.
Then the map $A\to\End_\C(M)$ is surjective.
\item Let $A,B$ complex algebras with units and let $R=A\otimes B$.
If $M,N$ are simple modules over  $A$ and $B$ resp., which are finite-dimensional over $\C$, then $M\otimes N$ is a simple $R$-module and every simple $R$-module, which is finite-dimensional over $\C$, is of this form with uniquely determined modules $M$ and $N$.
\end{enumerate}
\end{lemma}

\prf
Part (a) is a well-known result of Wedderburn.
See \cite{Lang} for a proof.

We show (b). According to  (a) it suffices to assume $A=\End_\C(M)$ and $B=\End_\C(N)$ for the first assertion.
The canonical map from $\End_\C(M)\otimes\End_\C(N)$ to $\End_\C(M\otimes N)$ is surjective, and so $M\otimes N$ is simple.

Let now $V$ be a given $\C$-finite-dimensional simple $A\otimes B$-module.
Then $V$ contains a simple $A$-module $M$ as $V$ is finite-dimensional.
Let $N=\Hom_A(M,V)$.
This vector space is a $B$-module in the obvious way.
Consider the map $\phi:M\otimes N\to V$ given by
$$
\phi(m\otimes \al)\=\al(m).
$$
Then $\phi\ne 0$ is an $A\otimes B$-homomorphism, hence surjective, as $V$ is simple.
We have to show that $\phi$ has zero kernel.
For this let $\al_1,\dots,\al_k$ be a basis of $N$ and $m_1,\dots,m_l$ a basis of $M$.
Let $c_{i,j}\in\C$ be given with $\phi\(\sum_{i,j}c_{i,j}m_i\otimes\al_j\)=0$.
We have to show that all coefficients $c_{i,j}$ are zero.
We have
$$
0\=\phi\(\sum_{i,j}c_{i,j}m_i\otimes\al_j\)\=\sum_{i,j}c_{i,j}\al_j(m_i)\=\sum_j\al_j\(\sum_ic_{i,j}m_i\).
$$
Let $P:M\to M$ be a projection onto a one-dimensional subspace, say $\C m_0$.
By part (a) there is $a\in A$ with $am=Pm$ for every $m\in M$.
We conclude
$$
0\= \sum_j\al_j\(P\(\sum_ic_{i,j}m_i\)\)\=\sum_j\la_j\al_j(m_0),
$$
where $P\(\sum_ic_{i,j}m_i\)\=\la_j m_0$.
For $a\in A$ arbitrary, we get
$
0\=\sum_j\la_j\al_j(am_0).
$
As $a$ runs through $A$, the vector $am_0$ runs through $M$, so
$
\sum_j\la_j\al_j\= 0.
$
Since the  $\al_j$ are linearly independent, all $\la_j$ are zero, hence all $\sum_ic_{i,j}m_i$ are zero and again by linear independence, we get $c_{i,j}=0$.
This proof also shows that all simple $A$-submodules of $V$ are isomorphic, whence the uniqueness.
\qed

We now prove Theorem  \ref{4.5.9}.
Let $\eta$ be a  $K\times L$-admissible representation of  $G\times H$.
Let $E$ be a finite subset of $\hat K$ and $F$ a finite subset of $\hat L$.
Then $V_\eta(E\times F)$ is a finite-dimensional irreducible $C_{E\times F}$-module.
By Lemma \ref{Lem7.5.22} the space  $V_\eta(E\times F)$ is an irreducible $C_E\otimes C_F$-module and by Lemma \ref{Lem7.5.23} the module $V_\eta(E\times F)$ is a tensor product of modules, which we write as $V_\pi(E)\otimes V_\tau(F)$.
The uniqueness of the tensor-factors give injective homomorphisms
$\ph_E^{E'}:V_\pi(E)\to V_\pi(E')$ if $E\subset E'$ and likewise for $F$.
By the uniqueness of the inner products in the Lemma of Schur, we can scale these homomorphisms in a way that they are isometric.
We define 
$$
\tilde V_\pi\df\lim_\to V_\pi(E).
$$
Then $\tilde V_\pi$ is a pre-Hilbert-space and we write
$V_\pi$ for its completion.
Analogously we construct $V_\tau$.
For every finite subset $E\subset \hat K$ the algebra  $C_E$ acts on $\tilde V_\pi$ by continuous operators.
These extend continuously to $V_\pi$ and we get *-representation of the Hecke-algebra $\CH_G$ and likewise for $V_\tau$.
By construction the isotypes of  $V_\pi$ are precisely the spaces  $V_\pi(E)$ for $E\in I_K$.
The isometrical maps $V_\pi(E)\otimes V_\tau(F)\hookrightarrow V_\eta$ induce an isometric
$\CH_G\otimes \CH_H$-homomorphism $\Phi:V_\pi\otimes V_\tau\to V_\eta$, which by irreducibility must be an isomorphism.
It only remains to install a unitary representation of the group $G\times H$ on $V_\pi\otimes V_\tau$, so at first a representation $\pi$ of $G$ on $V_\pi$.
Fix an arbitrary vector $w\in V_\tau$ with $\norm w=1$.
Then the map $V_\pi\to V_\pi\otimes w\tto\Phi V_\eta$ is an isometric embedding of $V_\pi$ into $V_\eta$, which commutes with the $\CH_G$-operation.
The  $G$-representation on $V_\eta$ then defines a unitary $G$-representation on $V_\pi$, which induces the $\CH_G$-representation.
We do the same with  $H$ and get irreduzible
unitary representations $\pi$ and $\tau$ of $G$ and $H$ such that $\Phi$ is a $G\times H$-isomorphism.
\qed

{\small Mathematisches Institut,
Auf der Morgenstelle 10,
72076 T\"ubingen,
Germany\\
\tt deitmar@uni-tuebingen.de}


\begin{thebibliography}{XXX}

\bibitem{HA2}
\bf Deitmar, Anton; Echterhoff, Siegfried:
\it Principles of harmonic analysis.
\rm Universitext. Springer, New York, 2009.

\bibitem{Flath}
\bf Flath, D.:
\it Decomposition of representations into tensor products.
\rm Automorphic forms, representations and $L$-functions (Proc. Sympos. Pure Math., Oregon State Univ., Corvallis, Ore., 1977), Part 1, pp. 179--183, Proc. Sympos. Pure Math., XXXIII, Amer. Math. Soc., Providence, R.I., 1979.

\bibitem{Lang}
\bf Lang, Serge:
\it Algebra.
\rm Revised third edition. Graduate Texts in Mathematics, 211. Springer-Verlag, New York, 2002

\end{thebibliography}
\end{document}